\documentclass{elsart}

\usepackage{amssymb,amsmath}

\newtheorem{theorem}{Theorem}[section]
\newtheorem{lemma}[theorem]{Lemma}

\newtheorem{proposition}[theorem]{Proposition}
\newtheorem{corollary}[theorem]{Corollary}

\newtheorem{question}[theorem]{Question}
\newtheorem{example}[theorem]{Example}

\journal{BULLETIN of the Australian Math. Soc.}
\begin{document}

\begin{frontmatter}

\title{Extension of continuous mappings and ${\rm H_1}$-retracts}

\author{Olena Karlova}

\address{Chernivtsi National University, Department of Mathematical Analysis,
Kotsjubyns'koho 2, Chernivtsi 58012, Ukraine\\mathan@ukr.net}

\begin{abstract}
We prove that any continuous mapping $f:E\to Y$ on a completely
metrizable subspace $E$ of a perfect paracompact space $X$ can be
extended to a Lebesgue class one mapping $g:X\to Y$ (i.e. for
every open set $V$ in $Y$ the preimage $g^{-1}(V)$ is an
$F_\sigma$-set in $X$) with values in an arbitrary topological
space $Y$.

\end{abstract}

\end{frontmatter}

\section{Introduction}

A mapping $f:X\to Y$ from a topological space $X$ to a topological
space $Y$ is called {\it a Lebesgue class $\alpha$ mapping} (or
{\it a mapping of the $\alpha$-th Lebesgue class}) if for every
closed set $F$ in $Y$ the set $f^{-1}(F)$ is of the multiplicative
class $\alpha$ in $X$. The family of all such mappings $f:X\to Y$
we denote by $H_\alpha(X,Y)$. Besides, we write $f\in H_1(X,Y)$ if
for every open set $V$ in $Y$ the preimage $g^{-1}(V)$ is an
$F_\sigma$-set in $X$.

Obviously, if $X$ or $Y$ is a perfect space then any continuous
mapping \mbox{$f:X\to Y$} is of the first Lebesgue class.

Classification of mappings  naturally leads to the problem on the
extension of mappings from a subset of a topological space to the
whole space with preserving of mapping's class or with it's
estimation. So, such classical results as Tietze Theorem
\cite[p.~116]{Eng} or Dugundji Theorem \cite{Dug} give the
possibility of extension of continuous mapping to continuous one.

Many mathematicians (F.~Hausdorff, W.~Sierpi\'{n}ski, M.~Alexits,
M.~Hahn, K.~Kuratowski) dealt in extension of real-valued
functions of some Lebesgue class.

K.~Kuratowski \cite{Ku2} proved that every mapping $f\in
H_{\alpha}(E,Y)$ on a subset $E$ of a metric space $X$ with values
in a complete metric separable space $Y$ can be extended to a
mapping $g:B\to Y$ of class $\alpha$ such that the set $B\supseteq
E$  is of the multiplicative class $\alpha+1$. Besides, if $E$ is
of the multiplicative class $\alpha>0$ then $f$ can be extended to
a Lebesgue mapping of class $\alpha$ on the whole space $X$.
Consequently, every mapping of the Lebesgue class $\alpha\ge 0$ on
a set $E\subseteq X$ can be extended to a mapping of the Lebesgue
class $\alpha+1$ on $X$. In particular, the following result
holds.

\begin{theorem}\label{Ku}

Let $X$ be a metric space, $Y$ be a complete metric separable
space and $E\subseteq X$. Then every continuous mapping $f:E\to Y$
can be extended to a mapping $g\in H_1(X,Y)$.

\end{theorem}

It follows from the above that the problem on the extension of
continuous function (Lebesgue class one function) to continuous
function (Lebesgue class one function) essentially differs from
the problems on the extension of functions preserving its class.
For example, if $X=\mathbb R$ and $E=\mathbb Q$ then not every
continuous function $f:E\to\mathbb R$ can be extended to a
continuous function defined on $X$; and it is easy to construct an
everywhere discontinuous function $f:E\to\mathbb R$ which is of
the first Lebesgue class and cannot be extended to a function of
the first Lebesgue class on $X$. On the other hand, Theorem
\ref{Ku} implies that every continuous function $f:E\to\mathbb R$
can be extended to a function $g:X\to\mathbb R$ of the first
Lebesgue class.

In connection with Theorem \ref{Ku} the following question arises.

\begin{question}\label{QKu}

Is it possible to omit the assumption of separability on space~$Y$
in Theorem \ref{Ku}?

\end{question}

R.~Hansell studied the problem of the extension of Lebesgue
mappings with non-separable metrizable ranges using the notion of
$\sigma$-discrete mapping introduced by A.~Stone \cite{Stone1}.

Recall that a family ${\mathcal A}$ of subsets of a topological
space $X$ is called {\it discrete} if for every point $x\in X$
there exists a neighborhood $U$ which intersects with at most one
set from ${\mathcal A}$.

A family ${\mathcal A}$ is called {\it $\sigma$-discrete} if it
can be written as a countable union of discrete families.

The family ${\mathcal B}$ of subsets of a topological space $X$ is
said to be {\it a base for a mapping} \mbox{$f:X\to Y$} if for
every open set $V$ in $Y$ there exists a subfamily ${\mathcal
B}_V\subseteq {\mathcal B}$ such that $f^{-1}(V)=\bigcup {\mathcal
B}_V$. If, moreover, the system ${\mathcal B}$ is
$\sigma$-discrete then it is called {\it a $\sigma$-discrete base
for $f$} and $f$ is called {\it a $\sigma$-discrete mapping}. The
family  of all $\sigma$-discrete mappings we denote by
$\Sigma(X,Y)$.

Evidently, every mapping with a second countable range space is
$\sigma$-discrete. Also it is easy to see that every continuous
mapping with metrizable domain or range is $\sigma$-discrete since
a metrizable space has a $\sigma$-discrete base \cite{Eng}.

The mentioned paper of R.~Hansell \cite{Ha2} contains the
following result.

\begin{theorem} {\rm \bf \cite[Theorem 9]{Ha2}\label{Ha}}
Let $X$ be a paracompact space, $Y$ a complete metric space,
$E\subseteq X$ and \mbox{$f:E\to Y$} a $\sigma$-discrete Lebesgue
mapping of class $\alpha$. Then $f$ can be extended to a Lebesgue
mapping $g:B\to Y$ of class $\alpha$ so that the set $B\supseteq
E$ is of multiplicative class $\alpha+1$.

\end{theorem}

The following question naturally arises.

\begin{question}\label{QHa}

Is it possible to replace the set $B$ in Theorem \ref{Ha} with the
whole space $X$?

\end{question}

The most recent, as to the best of our knowledge, result on the
extension of Lebesgue functions is due to O.~Kalenda and
J.~Spurn\'{y}.

\begin{theorem} {\rm\bf \cite[Theorem 29]{Ka}\label{Kalenda1}}
Let $E$ be a Lindel\"{o}f subspace of a completely regular space
$X$, $Y$ a complete metric separable space and

(i) $E$ be hereditarily Baire, or

(ii) $E$ be $G_\delta$ in $X$.

Then every mapping $f\in H_1(E,Y)$ can be extended to a mapping
$g\in H_1(X,Y)$.

\end{theorem}

At the same time it is interesting to study when can we extend
mappings with values in an arbitrary topological space.

In Section 2 we introduce and study the notion of $H_1$-retract
which is tightly connected with the problem on the extension of
continuous mappings to Lebesgue class one mappings with values in
an arbitrary topological space (analogously, as the notion of a
retract connected with extension of continuous mappings with
preserving of continuity).

Further in Section 3 we prove that every continuous mapping
$f:E\to Y$ on a completely metrizable subspace $E$ of a perfect
paracompact space $X$ with values in an arbitrary topological
space $Y$ can be extended to a Lebesgue class one mapping $g:X\to
Y$. This result implies a positive answer to Question~\ref{QKu}.
Besides, we give a negative answer to Question \ref{QHa}.

\section{$H_1$-retracts and their properties}

Let $X$ be a topological space and $E\subseteq X$. Recall
\cite{Bors} that a set $E$ is said to be {\it a retract} of $X$ if
there exists a continuous mapping $r:X\to E$ such that $r(x)=x$
for all $x\in E$. The mapping $r$ is called {\it a retraction} of
$X$ onto $E$. It is easy to see that a set $E\subseteq X$ is
retract of $X$ if and only if for any topological space $Y$ every
continuous mapping $f:E\to Y$ can be extended to a continuous
mapping $g:X\to Y$.

A subset $E$ of a topological space $X$ we call {\it an
$H_1$-retract} of $X$ if there exists a mapping \mbox{$r\in
H_1(X,E)$} such that $r(x)=x$ for all $x\in E$. The mapping $r$ we
call {\it an $H_1$-retraction} of $X$ onto $E$.

The following properties of $H_1$-retracts immediately follow from
the definition.

\begin{proposition}\label{prop1.1}

Let $X$ be a topological space. A set $E\subseteq X$ is
$H_1$-retract of $X$ if and only if for an arbitrary space $Y$
every continuous mapping $f:E\to Y$ can be extended to a Lebesgue
class one mapping $g:X\to Y$.

\end{proposition}

\begin{proposition}\label{propPerfect}

Let $E$ be an $H_1$-retract of a topological space $X$. Then $E$
is a perfect space.

\end{proposition}

A subset $A$ of a topological space $X$ is said to be {\it an
ambiguous set} if $A$ is simultaneously $F_\sigma$ and $G_\delta$
in $X$.

\begin{proposition}\label{prop1.2}

Let $X$ be a metrizable space and $E$ be an $H_1$-retract of $X$.
Then $E$ is $G_\delta$ in~$X$.

\end{proposition}

{\bf Proof.} Let $r:X\to E$ be an $H_1$-retraction of $X$ onto
$E$. It is easy to see that $E=\{x\in X:r(x)=x\}$.

Consider the diagonal $\Delta=\{(x',x'')\in X\times X: x'=x''\}$
of the space $X^2$ and the mapping $h:X\to X\times X$,
$h(x)=(r(x),x)$. Since $r\in H_1(X,E)\subseteq H_1(X,X)$ and the
mapping $g:X\to X$, $g(x)=x$ is continuous, according to
\cite[Theorem 1]{Ha1} the mapping $h:X\to X\times X$ is of the
first Lebesgue class. Since $\Delta$ is closed in $X\times X$, the
set $E=h^{-1}(\Delta)$ is $G_\delta$ in $X$.\hfill$\Box$

Point out that an $H_1$-retract may be in general even a
non-measurable set. Besides, the following example shows that the
assumption of metrizability of $X$ in the previous proposition is
essential.

\begin{example}
There exists a non-measurable $H_1$-retract $E$ of a perfect
separable linear ordered compact space $X$.
\end{example}

{\bf Proof.} Let $X=[0,1]\times \{0,1\}$ be endowed with the
lexicographic order, that is \mbox{$(x,i)<(y,j)$} if $x<y$ or
$x=y$ and $i<j$, $i,j\in \{0,1\}$. Remark that $X$ satisfies
necessary conditions (see \cite[P.~318]{Eng}).

Consider a set $E=\{(x,0):x\in [0,1]\}$. A mapping $r:X\to E$,
$r(x,i)=(x,0)$, is of the first Lebesgue class.

It remains to prove that $E$ is non-measurable.

For a set $A\subseteq X$ denote $A^+=\{x\in [0,1]:(x,1)\in A\}$
and $A^-=\{x\in [0,1]:(x,0)\in A\}$. It is not hard to prove that
for any open or closed set $A$ in $X$ we have $|A+\Delta
A^-|\le\aleph_o$. This implies that $|B^+\Delta B^-|\le \aleph_o$
for any measurable set $B$. But $E^+=\O$ and $E^-=[0,1]$. Hence,
$E$ is a non-measurable set. \hfill$\Box$

\begin{proposition}\label{prop1.4}

Let $X$ and $Y$ be topological spaces, $E$ be an ambiguous subset
of $X$ and $f:E\to Y$ be a Lebesgue class one mapping. Then there
exists a Lebesgue class one mapping $g:X\to Y$ such that $g|_E=f$.

\end{proposition}

\begin{corollary}\label{cor1.5}

Let $X$ be a topological space and $E$ be a perfect ambiguous
subset of $X$. Then $E$ is $H_1$-retract of $X$.
\end{corollary}

A subset $E$ of a topological space $X$ we call {\it a ${\rm
Coz}_\delta$-set} if there exists a sequence of continuous
functions $f_n:X\to [0,1]$ such that
$E=\bigcap\limits_{n=1}^\infty f_n^{-1}((0,1])$. The complement to
a ${\rm Coz}_\delta$-set we call a ${\rm Zer}_\sigma$-set. A set
which is simultaneously ${\rm Coz}_\delta$ and ${\rm Zer}_\sigma$
we call a {\it functionally ambiguous set}.

\begin{proposition}\label{Sequence}

Let $E_1,\dots,E_n$ be disjoint $H_1$-retracts of topological
space $X$ and $E_i$ is ${\rm Coz}_\delta$ in $X$ for every
$i\in\{1,\dots,n\}$. Then the union $E=\bigcup\limits_{i=1}^n E_i$
is an $H_1$-retract of $X$.

\end{proposition}

{\bf Proof.} First we obtain that for every finite family of
disjoint ${\rm Coz}_\delta$-sets $E_1,\dots,E_n$ there exist
disjoint functionally ambiguous sets $B_1,\dots,B_n$ such that
$E_i\subseteq B_i$ for every $i\in\{1,\dots,n\}$ and
$X=\bigcup\limits_{i=1}^n B_i$.

Let $n=2$ and $E_1$, $E_2$ be disjoint ${\rm Coz}_\delta$-sets.
Then the complements $E_i^c=X\setminus E_i$, $i=1,2$, are ${\rm
Zer}_\sigma$ and $E_1^c\cup E_2^c=X$. From  \cite[Lemma 3.2]{K1}
it follows that there exist functionally ambiguous sets $B_1$ and
$B_2$ such that $B_1^c\subseteq E_1^c$, $B_2^c\subseteq E_2^c$,
$B_1^c\cup B_2^c=X$ and $B_1^c\cap B_2^c=\O$. Then $E_1\subseteq
B_1$, $E_2\subseteq B_2$, $B_1\cap B_2=\O$ and $B_1\cup B_2=X$.

Let $n>2$ and the assumption holds when we have $n-1$ set. There
exist disjoint functionally ambiguous sets $\tilde{B}_1,\dots,
\tilde{B}_{n-1}$ such that $E_i\subseteq \tilde{B}_i$ if $1\le
i\le n-2$, $E_{n-1}\cup E_n\subseteq \tilde{B}_{n-1}$ and
$X=\bigcup\limits_{i=1}^{n-1}\tilde{B}_i$. Moreover, there exist
disjoint functionally ambiguous sets $C$ and $D$ such that
$E_{n-1}\subseteq C$, $E_n\subseteq D$ and $C\cup D=X$. Set
$B_i=\tilde{B}_i$ for $i=1,\dots,n-2$,
$B_{n-1}=\tilde{B}_{n-1}\cap C$ and $B_n=\tilde{B}_{n-1}\cap D$.

Let $r_i:X\to E_i$ be $H_1$-retractions, $1\le i\le n$. For every
$x\in X$ define $r(x)=r_i(x)$ if $x\in B_i$ for some
$i\in\{1,\dots,n\}$. Clearly, $r\in H_1(X,E)$ and $r(x)=x$ if
$x\in E$.\hfill$\Box$

\section{Extension of continuous mappings to the first class mappings
from completely metrizable subspaces}

In this section we prove the main results of this paper. All
topological spaces will be considered to be Hausdorff.

We say that a family ${\mathcal A}=(A_i:i\in I)$ of sets $A_i$
{\it refines} a family \mbox{${\mathcal B}=(B_j:j\in J)$} of sets
$B_j$ if for every $i\in I$ there exists $j\in J$ such that
$A_i\subseteq B_j$. We write this ${\mathcal A}\preceq {\mathcal
B}$.

\begin{lemma}\label{lemma1.7}

Let $X$ be a perfect paracompact space and ${\mathcal G}$ be a
locally finite cover of $X$ by ambiguous sets. Then there exists a
disjoint locally finite cover of $X$ by ambiguous sets which
refines ${\mathcal G}$.

\end{lemma}

{\bf Proof.} Without loss of generality we may assume that
\mbox{${\mathcal G}=\{G_\alpha:0\le \alpha<\beta\}$}, where
$\beta$ is some ordinal.

Denote $A_o=G_o$. For every $0<\alpha<\beta$ let
\mbox{$A_\alpha=G_\alpha\setminus \bigcup\limits_{\xi<\alpha}
G_\xi$}. According to Michael Theorem \cite[P.~430]{Eng}, the set
$\bigcup\limits_{\xi<\alpha} G_\xi$ is ambiguous as a locally
finite union of ambiguous sets. Then the set $A_\alpha$ is also
ambiguous. Clearly, the family ${\mathcal
A}=(A_\alpha:0\le\alpha<\beta)$ is to be found.\hfill$\Box$

The next theorem is the main result of our paper.

\begin{theorem}\label{t1.8}

Let $X$ be a perfect paracompact space and $E\subseteq X$ be a
completely metrizable subspace of $X$. Then $E$ is an
$H_1$-retract of $X$.

\end{theorem}

{\bf Proof.} Let $d$ be a metric on $E$  such that $(E, d)$ is a
complete metric space and $d$ induce the topology in $E$.

For every $n\in\mathbb N$ consider a cover ${\mathcal V}_n$ of the
set $E$ by open balls with radius $\frac{1}{2^{n+2}}$. For every
ball $V\in \bigcup\limits_{n=1}^\infty {\mathcal V}_n$ choose an
open set $U_V$ in $X$ so that $V=E\cap U_V$.

For every $n\ge 1$ let ${\mathcal G}_n=(U_{V_1}\cap \dots\cap
U_{V_n}: V_1\in {\mathcal V}_1, \dots,V_n\in {\mathcal V}_n)$ and
$G_n=\bigcup\limits_{G\in {\mathcal G}_n} G$.

Since $X$ is perfect, the set $G_n$ is an $F_\sigma$-set in $X$.
It follows from~\cite[P.~457]{Eng} that $G_n$ is paracompact
space. Then there exists a locally finite in $G_n$ cover
${\mathcal U}_n$ of $G_n$ by open sets in $G_n$, which refines
${\mathcal G}_n$. According to Lemma \ref{lemma1.7}, there exists
a disjoint locally finite in $G_n$ cover of $G_n$ by ambiguous
sets in $G_n$, which refines ${\mathcal U}_n$. Remove from this
cover those sets which do not intersect with $E$ and denote this
new system by ${\mathcal W}_n$. Let $P_n=\bigcup {\mathcal W}_n$.
Note that $P_n\subseteq G_n$. Fix an arbitrary set $W_n$ from
${\mathcal W}_n$ and denote \mbox{$W_n^o=W_n\cup (X\setminus
P_n)$}.

Index the elements of the system $\{{\mathcal W}_n^o\}\bigcup
({\mathcal W}_n\setminus \{{\mathcal W}_n\})$ and obtain the
family ${\mathcal X}_n=(X_{n,i}:i\in I_n)$.

Constructed in such a way sequence $({\mathcal X}_n)_{n=1}^\infty$
of families ${\mathcal X}_n$ satisfies the following properties:

(i) $X=\bigcup\limits_{i\in I} X_{n,i}$;

(ii) $X_{n,i}\cap X_{n,j}=\O,\quad i\ne j$;

(iii) $X_{n,i}\cap E\ne\O$ for all $i\in I_n$.

(iv) family $(X_{n,i}\cap G_n: i\in I_n)$ is locally finite in
$G_n$;

(v) $|\{i\in I_n:X_{n,i}\setminus G_n\ne\O\}|=1$;

(vi) ${\rm diam} (X_{n,i}\cap E)\le \frac{1}{2^{n+1}}$ for every
$i\in I_n$.

Since all the elements of the system ${\mathcal W}_n$ are
ambiguous sets in open subset $G_n$ of a perfect space $X$, all
the elements of system ${\mathcal W}_n$ are ambiguous sets in $X$.
Besides, since ${\mathcal W}_n$ is locally finite in $G_n$,
Michael Theorem \cite[P.~430]{Eng} implies that $P_n$ also is
ambiguous set in $X$. This implies that

(vii) $X_{n,i}$ is ambiguous in $X$ for all $i\in I_n$.

For every $n\in\mathbb N$ let
$$
E_{n,i}=X_{n,i}\cap E
$$
and for all $i_1\in I_1,\dots,i_n\in I_n$ let
$$
B_{i_1\dots i_n}=E_{1,i_1}\cap E_{2,i_2}\cap\dots\cap E_{n,i_n},
$$
$$
C_{i_1\dots i_n}=X_{1,i_1}\cap X_{2,i_2}\cap\dots\cap X_{n,i_n}.
$$
Then:

(1) $E=\bigcup\limits_{i_1\in I_1,\dots,i_n\in I_n}B_{i_1\dots
i_n}$ and $X=\bigcup\limits_{i_1\in I_1,\dots,i_n\in
I_n}C_{i_1\dots i_n}$ for every $n\in\mathbb N$;

(2) $B_{i_1\dots i_n}\cap B_{j_1\dots j_n}=\O$ and $C_{i_1\dots
i_n}\cap C_{j_1\dots j_n}=\O$ if $(i_1,\dots,i_n)\ne
(j_1,\dots,j_n)$;

(3) if $m\ge n$ and $C_{i_1\dots i_n}\cap C_{j_1\dots j_m}\ne\O$
then $i_1=j_1,\dots,i_n=j_n$;

(4) $C_{i_1\dots i_n}\cap E=B_{i_1\dots i_n}$ for every
$n\in\mathbb N$ and $i_1\in I_1$, $i_2\in I_2$,...,$i_n\in I_n$;

(5) $B_{i_1\dots i_n}$ is an ambiguous set in $E$ and $C_{i_1\dots
i_n}$ is an ambiguous set in  $X$ for all $n$ and $i_1\in I_1$,
$i_2\in I_2$,...,$i_n\in I_n$.

Moreover,

(6) for every $n\in\mathbb N$ and any set $I'\subseteq
I_1\times\dots\times I_n$ the set $\displaystyle
A=\bigcup\limits_{(i_1,\dots,i_n)\in I'}C_{i_1\dots i_n}$ is
ambiguous in $X$.

According to (iv), for every $k\in \mathbb N$ the family
$(X_{k,i}\cap G_k:i\in I_k)$ is locally finite in $G_k$, and (v)
implies that the family $(X_{k,i}\cap (X\setminus G_k):i\in I_k)$
is locally finite in $X\setminus G_k$. Therefore, taking into
consideration that the sequence $(G_n)_{n=1}^\infty$ decreases, we
obtain that for an arbitrary set
$$
D\in \{X\setminus G_1, G_1\setminus G_2,G_2\setminus
G_3,\dots,G_{n-1}\setminus G_n,G_n\}=\{D_o,\dots,D_n\}
$$
and for every $k\in\{1,2,\dots,n\}$ the family $(X_{k,i}\cap
D:i\in I_k)$ is locally finite in $D$. Then we have that the
family $(C_{i_1,\dots,i_n}:i_1\in I_1,\dots,i_n\in I_n)$ is also
locally finite in $D$. Hence, the family $(C_{i_1,\dots,i_n}\cap
D:(i_1,\dots,i_n)\in I')$ is locally finite in $D$. Further, since
all the sets $D_o,\dots,D_n$ are ambiguous in $X$ and (5) holds,
all the sets $A_k=\bigcup\limits_{(i_1,\dots,i_n)\in
I'}C_{i_1,\dots,i_n}\cap D_k$ are ambiguous in $X$ and
$A=\bigcup\limits_{k=0}^n A_k$ is ambiguous set in $X$.

For every  $n$ and $i\in I_n$ choose an arbitrary point
$y_{n_i}\in E_{n_i}$. For every $x\in E$ let
$\psi_n(x)=y_{n,i_n}\,\, \mbox{if} \,\, x\in B_{i_1\dots i_n}$.
Note that according to (1) and (2), mappings $\psi_n:E\to E$ are
correctly defined. Show that the sequence $(\psi_n)_{n=1}^\infty$
uniformly converges to the identical mapping \mbox{$\psi:E\to E$},
$\psi(x)=x$.

Fix $x\in E$ and $n\in\mathbb N$. Then there exist $i_1\in
I_1$,...,$i_n\in I_n$ such that $x\in B_{i_1\dots i_n}$. Then
$\psi_n(x)=y_{n,i_n}$. Since $B_{i_1\dots i_n}\subseteq
E_{n,i_n}$,  $x\in E_{n,i_n}$ and $y_{n,i_n}\in E_{n,i_n}$.
According to (vi), we have that ${\rm diam}
E_{n,i_n}\le\frac{1}{2^{n+1}}$. Then
$$
d(\psi(x),\psi_n(x))=d(x,y_{n,i_n})\le\frac{1}{2^{n+1}}.
$$

Note that
$$
d(\psi_m(x),\psi_n(x))\le
\frac{1}{2^{n+1}}+\frac{1}{2^{n+1}}=\frac{1}{2^n} \,\,\mbox{for
all}\,\, m\ge n\,\,\mbox{and}\,\, x\in E. \eqno(*)
$$

For every $n$ and multi-index $(i_1\dots i_n)\subseteq I_1\times
\dots\times I_n$ denote
$$
\ell(i_1\dots i_n)=\max\{1\le k\le n:B_{i_1\dots i_k}\ne\O\}.
$$

For all $n\in\mathbb N$ and $x\in X$ let $r_n(x)=y_{\ell(i_1\dots
i_n),i_{\ell(i_1\dots i_n)}}$ if $x\in C_{i_1\dots i_n}$.
Properties (1) and (2) imply that all the mappings $r_n:X\to E$
are correctly defined.

Prove that the sequence $(r_n)_{n=1}^\infty$ satisfies
inequality~$(*)$ for all $x\in X$.

Let $x_o\in X$ and $m\ge n$. Then there exist  $i_1\in
I_1$,...,$i_n\in I_n$ and \mbox{$j_1\in I_1$},...,\mbox{$j_m\in
I_m$} such that $x_o\in C_{i_1\dots i_n}\cap C_{j_1\dots j_m}$.
Property (3) implies $i_1=j_1,\dots,i_n=j_n$.

If $B_{i_1\dots i_n}\ne \O$ then $r_n(x_o)=y_{n,i_n}$. Let
$k=\ell(j_1\dots j_m)$. Then $r_m(x_o)=y_{k,j_k}$. Clearly, $k\ge
n$. Choose any point $x\in B_{j_1\dots j_k}$. Since $B_{j_1\dots
j_k}=B_{i_1\dots i_n j_{n+1}\dots j_k}\subseteq B_{i_1\dots i_n}$,
$\psi_n(x)=y_{n,i_n}$ and $\psi_k(x)=y_{k,j_k}$. Inequality  (*)
implies that
$$
d(r_n(x_o),r_m(x_o))=d(y_{n,i_n},y_{k,j_k})=d(\psi_n(x),\psi_k(x))<\frac{1}{2^n}.
$$

If $B_{i_1\dots i_n}=\O$ then $\ell(i_1\dots i_n)=\ell(j_1\dots
j_m)$. Now we have that $r_n(x_o)=r_m(x_o)$ and
$d(r_n(x_o),r_m(x_o))=0$. Hence, sequence $(r_n)_{n=1}^\infty$
satisfies (*) for all $x\in X$.

Since $X$ is a completely metrizable space, there exists a mapping
$r:X\to E$ such that the sequence $(r_n)_{n=1}^\infty$ uniformly
converges to $r$ on $X$. Besides, since $r_n|_E=\psi_n$ and
$\psi(x)=\lim\limits_{n\to\infty}\psi_n(x)$ for all $x\in E$, we
have that $r|_E=\psi$, that is $r(x)=x$ for every $x\in E$.

Since a uniform limit of the Lebesgue class one mappings is a
Lebesgue class one mapping \cite[P.~395]{Ku1}, it remains to prove
that $r_n\in H_1(X,E)$ for all $n\in\mathbb N$.

Since for any $n\in\mathbb N$ and $i_1,\dots,i_n\in
I_1\times\dots\times I_n$ such that $C_{i_1,\dots,i_n}\ne\O$ the
mapping $r_n|_{C_{i_1,\dots,i_n}}$ is constant, we have that for
an arbitrary set $B\subseteq E$
$$
r_n^{-1}(B)=\bigcup\limits_{(i_1,\dots,i_n)\in
I'}C_{i_1,\dots,i_n},
$$
where $I'=\{(i_1,\dots,i_n)\in I_1\times\dots\times
I_n:r_n(C_{i_1,\dots,i_n})\subseteq B\}$. Therefore, according to
(6), $r_n^{-1}(B)$ is an ambiguous set in $X$. In particular, all
the mappings $r_n$ are of the first Lebesgue class.\hfill$\Box$

\begin{corollary}\label{CorGdelta}

Let $X$ be a completely metrizable space and $E\subseteq X$. The
set $E$ is an  $H_1$-retract of $X$ if and only if $E$ is
$G_\delta$ in $X$.

\end{corollary}

{\bf Proof.} {\it Sufficiency.} It immediately follows from
\ref{prop1.2}.

{\it Necessity.} According to Aleksandrov-Hausdorff Theorem
\cite[p.~407]{Eng}, the space $E$ is completely metrizable. Hence,
Theorem \ref{t1.8} implies that  $E$ is an $H_1$-retract
of~$X$.\hfill$\Box$

The following corollary gives a positive answer to Question
\ref{QKu}.

\begin{corollary} \label{t1.9}

Let $X$ be a metrizable space, $Y$ be a completely metrizable
space, $A\subseteq X$ and $f:A\to Y$ be a continuous mapping. Then
there exists a Lebesgue class one mapping $g:X\to Y$ such that
$g|_A=f$.

\end{corollary}

{\bf Proof.} Denote by $\hat{X}$ the completion of $X$. According
to \cite[p.~405]{Eng}, there exists a $G_\delta$-subset $\hat{A}$
of $\hat{X}$ and a continuous mapping $h:\hat{A}\to Y$ such that
$A\subseteq \hat{A}$ and $h|_A=f$. According to Corollary
\ref{CorGdelta}, $\hat{A}$ is an $H_1$-retract of $\hat{X}$. Then
there exists a mapping $\hat{h}\in H_1(\hat{X},Y)$ such that
$\hat{h}|_{\hat{A}}=h$.

Let $g=\hat{h}|_X$. Then $g:X\to Y$ is the desired extension of
$f$.\hfill$\Box$

Since every completely metrizable separable space is hereditarily
Baire and Lindel\"{o}f,  the result of O.~Kalenfa and
J.~Spurn\'{y} implies the following fact.

\begin{theorem}\label{Kalenda}

Let $E$ be a completely metrizable separable subspace of a
completely regular space $X$. Then $E$ is an $H_1$-retract of $X$.

\end{theorem}

At first sight this theorem gives a solution to the problem on the
extension of continuous mapping to a mapping of the first Lebesgue
class with values in an arbitrary (not necessary separable)
topological space, analogously as in Theorem~\ref{t1.8}. But since
a continuous image of a separable space $E$ is also separable, in
fact, separability of $Y$ is present here imperceptibly, and we
cannot obtain Corollary \ref{t1.9} from Theorem \ref{Kalenda}.

The following example shows that the assumption that $X$ is
perfect in Theorem  \ref{t1.8} and the assumption that $E$ is
separable in Theorem \ref{Kalenda} cannot be omitted. Besides,
this example gives a negative answer to Question \ref{QHa}.

\begin{example}

There exist a completely metrizable subspace $E$ of a compact
space $X$ and a continuous function $f:E\to [0,1]$ which cannot be
extended to a Lebesgue class one function on $X$.

\end{example}

{\bf Proof.} Let $E$ be an uncountable discrete space and
$X=\alpha E=E\cup\{\infty\}$ be the Aleksandrov compactification
of $E$.

Choose two uncountable disjoint subsets $E_1$ and $E_2$ of $E$ so
that $E=E_1\sqcup E_2$ and consider the function
$$
f(x)=\left\{\begin{array}{ll}
  1, & \mbox{if}\quad x\in E_1, \\
  0, & \mbox{if}\quad x\in E_2. \\
\end{array}
\right.
$$
The function $f:E\to [0,1]$ is continuous and hence a
$\sigma$-discrete function of the first Lebesgue class.

Remark that for every continuous function  (and for every Baire
one function, i.e. a pointwise limit of continuous functions)
$g:X\to [0,1]$ there exists at most countable set  $X_o\subseteq
X$ such that $g(x)=g(\infty)$ for all $x\in X\setminus X_o$. It
follows that a function $f$ cannot be extended to a Baire one
function $g:X\to [0,1]$, provided $E_1$ and $E_2$ are uncountable
sets.

According to  \cite[Theorem 3.7]{Ve}, the class $H_1(X,[0,1])$
coincides with the class of all Baire one functions $g:X\to
[0,1]$. Therefore, the function $f$ cannot be extended to a
Lebesgue class one function on $X$.\hfill$\Box$

Author would like to thank to professor V.~Mykhaylyuk for helpful
suggestions and comments.

\end{document}